\documentclass[11pt]{article}
\usepackage{amsfonts,amsmath,latexsym,verbatim,amscd,amsthm}
\newcommand{\be}{\begin{eqnarray}}     
\newcommand{\ee}{\end{eqnarray}}

\newcommand{\RR}{\mathbb{R}}

\newcommand{\ZZ}{\mathbb{Z}}

\newcommand{\e}{\epsilon}
\newcommand{\del}{\partial}

\newcommand{\calC}{{\mathcal C}}
\newcommand{\calE}{{\mathcal E}}

\newcommand{\calO}{{\mathcal O}}

\newcommand{\olR}{\overline{R}}
\newcommand{\hc}{\mathrm{cusp}}
\newcommand{\eucl}{\mathrm{euc}}

\title{Ricci flow on surfaces with cusps}
\author{Lizhen Ji \thanks{lji@umich.edu, University of Michigan, Ann Arbor, MI 48105}
\and Rafe Mazzeo \thanks{mazzeo@math.stanford.edu, Stanford University, Stanford, CA 94305}
\and Natasa Sesum \thanks{sesum@math.columbia.edu, Columbia University, New York, NY 10027}}

\date{} 
\theoremstyle{plain} 
\newtheorem{dummy}{Dummy}

\theoremstyle{definition}

\theoremstyle{plain}
\newtheorem{corollary}[dummy]{Corollary}
\newtheorem{remark}[dummy]{Remark}
\newtheorem{lemma}[dummy]{Lemma}
\newtheorem{theorem}[dummy]{Theorem}
\newtheorem{proposition}[dummy]{Proposition}

\begin{document}
\maketitle

\begin{abstract}
We consider the normalized Ricci flow $\del_t g = (\rho - R)g$ with initial condition a complete 
metric $g_0$ on an open surface $M$ where $M$ is conformal to a punctured compact Riemann surface 
and $g_0$ has ends which are asymptotic to hyperbolic cusps.  We prove that when $\chi(M) < 0$ and
$\rho < 0$, the flow $g(t)$ converges exponentially to the unique complete metric of constant
Gauss curvature $\rho$ in the conformal class. 
\end{abstract}

\section{Introduction}
The general uniformization theorem states that if $M$ is an arbitrary Riemann surface, i.e.\ 
a two-dimensional manifold (either open or closed) with a specified conformal structure, then up to 
a constant multiple there is a unique metric in the given conformal class which is complete and has 
constant curvature. Our interest here is with the case when $M$ is open, and in this setting, 
in all but a few cases, this `uniformizing metric' has constant curvature $-1$. In particular, 
this is the case if $M$ is conformally equivalent to the complement 
of $\ell$ points $\{p_1, \ldots, p_\ell\}$ in a compact Riemann surface $\overline{M}$, so long as
$\chi(M) = \chi(\overline{M}) - \ell < 0$. (In other words, we exclude only the sphere 
minus one or two points, in which case the uniformizing metric is flat.)  We henceforth assume that 
this is the case without further comment.

The behaviour of Ricci flow on compact surfaces is mostly understood \cite{Ch}, \cite{Ha}, 
but relatively little work has been done in the noncompact complete setting. The case which 
is somehow closest to the compact one is when $(M,g_0)$ has finite area and finite total curvature, 
and in particular has ends which are asymptotic to hyperbolic cusps. The uniformizing metric 
$g_\infty$ in this case is asymptotically equivalent to $g_0$ (hence in particular, quasi-isometric
to it); more specifically, there exists a unique complete metric $\overline{g} = \overline{u}\, g_0$
with constant negative  Gauss curvature, and with $\overline{u}$ converging to a constant at infinity
on each end. Thus it seems reasonable that the (normalized) Ricci flow starting at $g_0$ should converge 
asymptotically to this uniformizing metric, and this is exactly what we show.  We allow the initial metric to have different limiting curvatures
on each end; these curvatures even out in the limit and the flow converges to a metric of constant
curvature (where the constant is determined by the normalizing constant in the flow equation).

To state the main result more carefully, we recall that a standard formula for the exact hyperbolic 
cusp metric (with $K \equiv -1$) on a punctured disk $|z| < 1$ is $g_c = (|z|\log |z| )^{-2}|dz|^2$ 
(the boundary at $|z| = 1$ is at finite distance). We prefer to use two other standard forms:
\begin{equation}
g_c\quad  = \ 
\left\{
\begin{array}{rcl}
&(|x|\log |x|)^{-2}\, |dx|^2, \qquad &x \in \RR^2 \setminus B_1 \qquad \mbox{or}\\
&ds^2 + e^{-2s}d\theta^2,  \qquad &s \geq 0,\ \theta \in S^1.
\end{array}
\right.
\label{eq:cm}
\end{equation}
Let us say that a metric $g$ is asymptotic to a multiple of $g_c$ if on each end $E \cong [s_0,\infty) \times S^1$, 
\begin{equation}
g = u \, g_c, \qquad u = \frac{2}{\sigma} + v \qquad \mbox{where}\quad \sigma \in \RR^+\ \mbox{and}\ |v| \leq C s^{-\mu}
\label{eq:acm}
\end{equation}
for some $\mu > 0$. Note that with this convention, $R \to -\sigma$ on $E$. The constant $\sigma$ can 
be different on different ends $E$ of $M$. (We also assume that 
$u$ is continuous, and that its derivatives up to order $2+\alpha$ satisfy certain estimates which will be 
made precise below.) 

The normalized Ricci flow for a family of metrics $g(t)$ on $M$ is the parabolic differential equation
\begin{equation}
\frac{\del\,}{\del t}g(t) = (\rho - R)g; \quad g(0) = g_0,
\label{eq:nRf}
\end{equation}
where $R$ is the scalar curvature (for surfaces, twice the Gauss curvature) and $\rho$ is a negative constant.
One good choice for this constant is
\begin{equation}
\rho = \olR_0 = \frac{\int_M R_{g_0}\, dA_0}{\int_M dA_0},
\label{eq:nRfc}
\end{equation}
in which case the family of solution metrics $g(t)$ has constant area, as we show below,
but it will be convenient to consider the equation for other values of $\rho$ too. Throughout
this paper we shall use the overline to denote an average value, as here. Note that this
flow preserves the conformal class. There is a general short-term existence result due to Shi \cite{Sh}
for the Ricci flow on complete manifolds of any dimension; however, his result does
not guarantee that the quasi-isometry class is preserved, so this will be our first result:

\begin{theorem}
Let $(M,g_0)$ be a complete surface with finite area such that $g_0$ is asymptotic
to a multiple of a hyperbolic cusp metric on each end of $M$. Then (\ref{eq:nRf}) has a unique
solution, which exists on some time-interval $0 < t < \e$, and such that for each 
$t$ in this range, the solution metric $g(t)$ is also asymptotic to a scaled hyperbolic 
cusp metric on each end.
\end{theorem}
\begin{remark}
In \S 2 we introduce weighted H\"older spaces which measure the rate of decay of the
conformal factor to a constant; this gives precise meaning to the phrase ``asymptotic to a 
scaled hyperbolic cusp metric.'' 
\end{remark}

The main result of this paper is 
\begin{theorem}
Let $(M,g_0)$ be as above and suppose that $\chi(M) < 0$ and $\rho < 0$. Then the solution $g(t)$ 
to the normalized Ricci flow equation (\ref{eq:nRf}) with initial condition $g_0$ exists for all 
$t > 0$ and converges exponentially to the unique complete metric with constant scalar curvature $\rho$ 
in this conformal class.
\end{theorem}

This implies the uniformization theorem for this class of surfaces:
\begin{corollary}
If $\chi(M) < 0$ and $(M,g_0)$ is any surface with a complete metric where on each end $g_0$ is asymptotic to 
a hyperbolic cusp, then there is a unique complete metric $g_\infty$ in the same conformal class as $g_0$ which has
constant negative curvature $-1$. 
\end{corollary}

We note that if an open Riemann surface $M$ is conformally equivalent to the complement of 
finitely many points of a compact Riemann surface, then it admits a conformally related 
complete metric which is exactly hyperbolic on every cusp neighborhood. Therefore, 
if $\chi(M)<0$, then $M$ admits a unique complete hyperbolic metric in its conformal class.
There are, of course, many different proofs of this classical theorem, including quite a few which
rely primarily on PDE techniques. For proofs in the compact case using Ricci flow see \cite{Ha}, 
\cite{Ch} and \cite{CLT}; \cite{Chru} has a proof using a fourth order flow. The paper \cite{MT} 
contains a PDE proof of the general uniformization theorem on arbitrary open Riemann surfaces.
For more about uniformization, see \cite{Ab}. Given all of this, perhaps the best emphasis 
for our result here is that it proves the stability for Ricci flow around complete finite area 
hyperbolic metrics: the initial metric can be quite far from hyperbolic on any compact set, so long as 
it is asymptotically hyperbolic at infinity. 

It would be very satisfying to analyze the Ricci flow on complete open surfaces $(M,g_0)$
assuming only that $(M,g_0)$ has finite total curvature, in which case a famous theorem
due to Huber states that it is conformal to a punctured compact Riemann surface. Our result 
is a step in this direction, and handles some part of the case where $(M,g_0)$ also has 
finite total area. 

We proceed as follows: we discuss some generalities about Ricci flow on open surfaces in the 
next section, then go on to prove the short-time existence for the flow in the class of metrics 
asymptotic to scaled hyperbolic cusps. The proof of long-time existence and convergence uses ideas 
from \cite{Ha}, in particular the existence of a ``potential function'', i.e.\ a solution $f_0$ 
to the equation $\Delta_{g_0}f_0 = R_0 - \overline{R}_0$ with certain properties, which is the subject of \S 3. 
The main result is proved in \S 4. 

\bigskip

\noindent {\bf Acknowledgements:} The first author is partially supported by
the NSF grant DMS-0604878 and would like to thank S.Wolpert for helpful correspondence.
The second author is partially supported by the NSF grant DMS-0505709. The third author 
wishes to thank Albert Chau, Peter Li, John Lott and Jiaping Wang for many useful discusions
regarding the problem, and is partially supported by the NSF grant DMS-0604657; she also 
acknowledges the support and hospitality at MSRI during the period in which this work was
initially carried out. An anonymous referee provided useful comments on an earlier (and
very much different) version of this paper, and Pierre Albin and Frederic Rochon made further 
useful comments on a much more recent version.

\section{Short-time existence}
Let $g(t)$ be any family of metrics on a surface $M$ which is $\calC^1$ in $t$ and 
which remains in the same conformal class, so we can write $g(t) = u(t)g_0$, $u(0) = 1$. 
Suppose that $\del_t g = (\rho -R)g$, where $\rho$ is any fixed constant and $R$ is the
scalar curvature of $g$. Using the usual transformation rule for the scalar curvature,
\begin{equation}
\Delta_0  \log u + Ru = R_0,
\label{eq:transfsccurv}
\end{equation}
the flow equation can be written as the scalar equation
\begin{equation}
\frac{\del u}{\del t} = \Delta_0 \log u + \rho u - R_0.
\label{eq:tRf}
\end{equation}
Here and for the rest of the paper, a subscript $0$ indicates that the relevant
quantity is associated to the metric $g_0$; a subscript $t$ will sometimes be added for
emphasis that a quantity is associated with $g(t)$. 

Differentiating (\ref{eq:transfsccurv}) with respect to $t$ leads to the equation
\begin{equation}
\frac{\del R}{\del t} = \Delta R + R(R-\rho),
\label{eq:eqforR}
\end{equation}
since $\Delta = u^{-1}\Delta_0$. Another useful fact which follows directly from (\ref{eq:transfsccurv})
and (\ref{eq:tRf}) is that
\begin{equation}
\frac{\del\,}{\del t}\, dA_t = (\rho - R)\, dA_t.
\label{eq:tdA}
\end{equation}
These two formul\ae\  already indicate that if $g(t)$ converges to a metric $g_\infty$ as $t \to \infty$, 
then the scalar curvature of $g_\infty$ is equal to $\rho$ and (so long as the convergences of
$R(\cdot,t)$ to $\rho$ is sufficiently quick), the area of $(M,g(t))$ converges to a nonzero limit too.

The next point is that for consideration of the short-time problem, the value of the constant $\rho$ 
is not important. Indeed, there is a simple transformation between solutions corresponding to different 
values of $\rho$. It is enough to check that solutions of the equation with $\rho=0$ can be transformed 
to solutions of the equation with any other $\rho \neq 0$. Suppose that $\del_t g = -R g$, and define
\[
\tilde{g}(\tau) = e^{\rho \tau}g, \qquad \mbox{where} \quad \tau = -\frac{1}{\rho} \log (1-\rho t).
\]
A straightforward computation yields $\del_\tau \tilde{g} = (\rho - \tilde{R})\tilde{g}$, as desired. 

Because of this, it is enough to prove short-time existence for the equation with the parameter $\rho$
set equal to $0$. As already noted, afterwards we shall want to set $\rho$ equal to any negative constant.
The precise value turns out not to matter too much, although the choice 
\begin{equation}
\rho = \olR_0 = \frac{4\pi \chi(M)}{\mbox{Area}\,(M,g_0)}
\label{eq:correctr}
\end{equation}
is common in the literature since the flow is then area-preserving. 

We now define the function spaces which will be used below. First, let $\calC^{k,\alpha}_\hc(M)$ denote 
the ordinary H\"older space, defined with respect to any fixed metric $g$ on $M$ which equals the model 
hyperbolic cusp metric $g_c$ on each end. In other words, derivatives and difference quotients
are defined with respect to the metric $g$. To be specific, the H\"older seminorm is defined by the
usual formula, but taking the supremum over all points $x,x' \in M$ with $\mbox{dist}\,(x,x') \leq 1$.
The fact that the injectivity radius tends to zero on each end $E$ does not cause a problem, since we
can pass to the universal cover of $E$ and compute the H\"older seminorms over unit balls there. 
Next, fixing coordinates $(s_j,\theta_j)$ on each end $E_j$ as in (\ref{eq:cm}) and a weight parameter 
$\mu > 0$, define the weighted H\"older spaces
\[
s^{-\mu} \calC^{k,\alpha}_\hc(M) = \{u \in \calC^{k,\alpha}_{\mathrm{loc}}(M): \left. u \right|_{E_j} = 
s_j^{-\mu}v_j, v_j \in \calC^{k,\alpha}_\hc(M)
\ \mbox{for}\ s_j \geq \bar{s} \}.
\]
Finally, if $M$ has $\ell$ ends, choose a partition of unity $\{\chi_0, \chi_1, \ldots, 
\chi_\ell\}$ where $\chi_0$ is compactly supported and equal to $1$ on some large compact set 
and $\chi_j$, $j \geq 1$, equals $1$ sufficiently far out on the end $E_j$, and define the map
\begin{equation}
\RR^\ell_+ \in \lambda = (\lambda_1, \ldots, \lambda_\ell) \longmapsto 
\kappa(\lambda) = \chi_0 + \sum_{j=1}^\ell \chi_j \lambda_j.
\label{eq:kappa}
\end{equation}
This family of functions will be used repeatedly below to rescale the metric on each end so as to be able to 
vary the limiting curvatures independently. 

We now state and prove our first result:
\begin{proposition}
Let $(M,g_0)$ be a surface with scaled asymptotically cusp ends in the sense that on each end,
\[
\left. g_0 \right|_{E_j} = w_{0,j} g_c, \qquad \mbox{where} \qquad
w_{0,j} := \frac{2}{\sigma_j}(1 + \tilde{w}_{0,j}),
\]
with $\chi_j \tilde{w}_{0,j} \in s^{-\mu}\calC^{2,\alpha}_\hc(M)$ for some $\mu > 0$ (so $R_0 \to -\sigma_j$ on that end). 
Then there is a unique solution $g(t) = u(t)g_0$ to the unnormalized Ricci flow (\ref{eq:tRf}) 
with $\rho =0$, with initial condition $u(0) = 1$ on some short time interval $0 \leq t < \e$. Moreover, 
for each $t$ in this interval, $u(t) = \kappa(\lambda(t)) + v(t)$ with $v(t) \in s^{-\mu}\calC^{2,\alpha}_\hc(M)$, 
$\lambda_j(0) = 1$ for all $j$ and $v, \lambda_1, \ldots, \lambda_\ell \in \calC^1$ as functions of $t$.
\label{pr:ste}
\end{proposition}

We only sketch the proof briefly, since it is elementary in this two-dimensional setting, and since 
short-time existence for this equation is already proved in \cite{Sh}, \cite{Sh1}.  The main point of 
Proposition \ref{pr:ste}, however, is that the asymptotic type of the metric is preserved under the flow, 
which is not proved by Shi.

It is convenient to use the change of variable $u = e^\phi$ (and recall again that $\rho=0$ here) 
to rewrite (\ref{eq:tRf}) as
\begin{equation}
\frac{\del \phi}{\del t} = e^{-\phi}\Delta_0 \phi - R_0 e^{-\phi}.
\label{eq:chvar}
\end{equation}
A solution to this quasilinear equation is obtained by a standard contraction mapping argument in the space 
$\calC^1([0,\e); \RR^\ell_+ \times s^{-\mu}\calC^{2,\alpha}_\hc(M))$, where $\lambda \in \RR^\ell_+$ is 
identified with the smooth function $\kappa(\lambda)$ which is constant on each end. To carry this out, one 
must understand the mapping properties for the linearization of (\ref{eq:chvar}) at $\phi = 0$ ($u = 1$), i.e.\ for 
the linear heat equation
\[
\del_t \psi =  \Delta_0 \psi + R_0 \psi.
\]
To this end, first note that by (\ref{eq:transfsccurv}) (and using the notation from the statement of this Proposition),
\[
\sigma_j + R_0 =  - \frac{\sigma_j}{2}\, \frac{1}{1 + \tilde{w}_{0,j}}\, \Delta_c \log (1 + \tilde{w}_{0,j}) 
\in s^{-\mu}\calC^{0,\alpha}_\hc(M).
\]
The main issue is to show that the linear flow acts on functions of the form $\kappa(\lambda(t)) + v(t, z)$
with $v$ satisfying the polynomial decay condition on $M$, provided $\lambda(t)$ is chosen correctly.
Indeed, from (\ref{eq:tRf}) we must set $\lambda_j(t) = 1 - \sigma_j t$ for each $j$ to make
the asymptotically constant part of this equation hold on each $E_j$. (Note that the decomposition of the conformal
factor into the sum of a constant and decaying term on each end makes sense for small times, but does
not hold uniformly as $t$ gets larger.) To control the other term, we use that the function $s^{-\mu}$ 
is a supersolution to the equation for $s \geq s_0$ for some sufficiently large $s_0$. Indeed, 
in this cusp region, 
\[
\Delta_0 + R_0 = w_0^{-1}\Delta_c + R_0 = w_0^{-1}\left(\del_s^2 - \del_s + e^{-2s}\del_\theta^2\right) + R_0
\]
so
\[
(\Delta_0 + R_0)s^{-\mu} = w_0^{-1}(\mu s^{-\mu-1} + \mu(\mu+1) s^{-\mu-2}) + R_0 s^{-\mu} < 0
\]
whenever $R_0 < -1/2$ and $s_0$ is sufficiently large. From this we define, for $A > 0$,
\[
\psi = C\, e^{At}\min\{1, (s/s_0)^{-\mu}\}
\]
on each end and $\psi = C$ away from the ends where $C$ is some large constant; this satisfies
\[
\left(\del_t - (\Delta_0 + R_0)\right) \psi \geq C'\, e^{At}\, \min\{1, (s/s_0)^{-\mu}\}
\]
for some $C' > 0$ (which depends only on $A, s_0$ and $R_0$), uniformly for $0 \leq t \leq t_0$, 
hence is a global supersolution for the parabolic equation. 

The rest of the proof uses the solution operator for this linear problem, Duhamel's principle 
and a contraction mapping argument to produce a unique solution of (\ref{eq:chvar}) in some short
time interval $[0,t_0]$. Since the equation is quasilinear, one needs to use the explicit structure
of the nonlinearity to make sure that one really does obtain a contraction. We leave the details to the reader. 

We conclude this section with the following computation. Let $(M,g_0)$ be a complete metric with scaled 
asymptotically hyperbolic cusp ends (with $R_0 \to -\sigma_j$ on the end $E_j$), and suppose that 
$\chi(M) < 0$. Let $g(t)$ be any solution of the 
normalized Ricci flow equation with $\rho < 0$ with initial condition $g_0$. Write 
$g(t) = (\kappa(\lambda(t)) + v) (\kappa(2/\sigma) + w_0) g_c$, where $2/\sigma = (2/\sigma_1, \ldots,
2/\sigma_\ell)$ and $v, w_0$ tend to $0$ like $s_j^{-\mu}$ along with their derivatives on each end. 
Then the Gauss-Bonnet theorem holds, i.e.\ 
\[
\int_M R\, dA = 4\pi \chi(M).
\]
To verify this, apply the Gauss-Bonnet theorem for surfaces with boundary to the 
surface $M$ with each cusp truncated at $s = s_0$ and note that the boundary terms 
in this formula decay as $s_0 \to \infty$. 

Next, let $A(t)$ denote the area of $(M,g(t))$; then by what we have just shown, 
\[
A'(t) = \int_M (\rho - R(\cdot,t))\, dA_t = \rho A(t) - 4 \pi \chi(M).
\]
Since both $\chi(M)$ and $\rho$ are negative, the value $A_\infty = 4\pi \chi(M)/\rho$ 
is a stable fixed point for this ODE, and hence represents the limiting value of $A(t)$
if the flow exists for all $t \geq 0$. If $\rho$ is given by the precise value (\ref{eq:correctr}),
then this ODE can be written as
\[
A'(t) = -4\pi \chi(M)\, \left( 1 - \frac{A(t)}{A(0)} \right),
\]
hence $A(t) \equiv A(0)$ for all $t$. Note that the behaviour of this ODE is quite different 
when $\chi(M)$ is nonnegative. 

To simplify matters below, we assume henceforth that $\rho = \olR_0 = 4\pi \chi(M)/A_0$, the average value 
of the scalar curvature of $g_0$, so the flow equation becomes 
\begin{equation}
\del_t g = (\olR -R) g; 
\label{equation-RF}
\end{equation}
clearly $\olR_t = \rho$ for any $t \geq 0$. Note that if we prove long-time existence and convergence of the
equation for this particular value of $\rho$, then using the transformation described at the beginning of this 
section and uniqueness for the solution, we obtain long-time existence and convergence for the equation with any 
$\rho < 0$; in other words, ultimately the value of $\rho$ is unimportant so long as it is negative.

\section{The potential function}
A clever idea, going back to Hamilton \cite{Ha}, to prove long-time existence of the flow
relies on the existence of a `potential function', which by definition is a function satisfying
\begin{equation}
\Delta f = R - \olR
\label{eq:potfcn}
\end{equation}
at time $t=0$, and which has vanishing integral and bounded gradient. Here $\olR$ is the average 
of the scalar curvature $R$ over $M$. The goal of this section is to find this function $f$.

\begin{proposition}
Let $(M,g)$ be metric on the surface $M$ with scaled asymptotically cusp ends, so $g = u g_c$ with 
$u = \kappa(2/\sigma) + v$, $v \in s^{-\mu}\calC^{2,\alpha}_\hc(M)$ for some fixed $\mu > 1$. Then there is a unique 
function $f$ of the form
\[
f = \sum_{j=1}^\ell \chi_j \left( c_j \log \log r_j + \beta_j\right) + 
\tilde{f}, \qquad \tilde{f} \in \calC^{2,\alpha}_\hc(M),
\]
for some constants $c_j, \beta_j$, which satisfies (\ref{eq:potfcn}), and in addition 
\[
\int_M f\, dA = 0, \qquad \mbox{and}\qquad \sup_M |\nabla f| < \infty.
\]
\label{pr:potential}
\end{proposition}

There are two obvious ways to attempt to find a solution to \eqref{eq:potfcn}: either via the standard $L^2$
existence theory (which by self-adjointness requires that the right hand side has integral equal to $0$),
or else by finding the solution as a limit of solutions on an expanding set of compact truncations of $M$. 
Unfortunately, neither approach works directly: the solution obtained by the first approach 
(though it does turn out to be the correct one) does not automatically have bounded gradient; 
on the other hand, in order to obtain the correct solution using the second approach, one must 
choose the boundary conditions on the various receding boundary components of these truncations carefully.
This is more or less equivalent to what takes place in the proof below. 

The proof has several steps. We shall write $g$ as a conformal factor times an asymptotically
Euclidean metric and use the more standard and elementary mapping properties for the corresponding 
Laplacian; we also first assume that $g$ decays very quickly to some scaling of $g_c$ at each end, 
and at the end deduce the result with the relatively slower decay stated here. 

Assuming as usual that $R \to -\sigma_j$ on $E_j$, we can write $g = \kappa(2/\sigma) A g_c$ where 
$g_c = (|x| \log |x|)^{-2}|dx|^2$ on $\RR^2\setminus B_c$ for $c \gg 0$. We denote by $g_e$ 
a metric conformal to $g$ which equals the Euclidean metric $|dx|^2$ on each end, for $|x| \geq c$, and 
which agrees with $g$ on some large compact set in $M$. Thus $g = \kappa(2/\sigma) A u_c g_e$ where 
$u_c = (|x| \log |x|)^{-2}$ on each (Euclidean) end. Define weighted H\"older spaces 
$r^{-\nu}\calC^{k,\alpha}_\eucl (M)$ with respect to the radial function $r = |x|$ on the 
Euclidean end and where derivatives are taken relative to the metric $g_e$ 
(rather than $g_c$ as in \S 2). If $M$ has more than one end, we sometimes let $r_j$ denote the radial 
function on the $j^{\mathrm{th}}$ end. We assume that 
\[
A - 1  := B \in r^{-\nu} \calC^{2,\alpha}_\eucl (M)\quad \mbox{for some} \quad \nu \in (0,1).
\]
Since $r = \exp(e^s)$ where $s$ is the geodesic distance for $g_c$, this is a {\it much} stronger 
decay condition than we had assumed earlier; it will be removed at the end of the proof.

As before, using (\ref{eq:transfsccurv}), since $\Delta_e 1 = \Delta_e \log r = 0$ and 
$\Delta_e \log \log r = -(r\log r)^{-2}$ and $R_e \equiv 0$ on $E_j$, we compute that 
\begin{multline}
\sigma_j + R = \sigma_j - \frac{\sigma_j}{2A} (r \log r)^2 \Delta_e (\log A - 2 \log \log r)  \\ 
= \sigma_j \frac{B}{1+B} + \frac{\sigma_j(r\log r)^2}{2(1+B)} \Delta_e \log(1+B).
\label{eq:RonEj}
\end{multline}
Therefore, $|\sigma_j + R| \leq C r_j^{-\nu}(\log r_j)^2$ on $E_j$, with the corresponding bound for its 
H\"older quotient. Now rewrite (\ref{eq:potfcn}) as
\[
\Delta_e f = q := \kappa(2/\sigma) A (r\log r)^{-2}(R - \olR). 
\]
By (\ref{eq:RonEj}), 
\[
\left. q \right|_{E_j} = \frac{2}{\sigma_j} A (r\log r)^{-2}\left( (R + \sigma_j) - (\sigma_j + \olR) \right)
= - c_j (r\log r)^{-2} + \Delta_e \log (1+B) + q_j,
\]
where $q_j \in r^{-\nu-2}\calC^{2,\alpha}_\eucl(E_j)$ and $c_j \in \RR$. Thus if we set
\[
f = \sum_{j=1}^\ell c_j \chi_j \log \log r_j + \log (1+B) + h
\]
then we have reduced (\ref{eq:potfcn}) to the problem of finding a solution $h$ to $\Delta_e h = \tilde{q}$,
where $\tilde{q} := q + \sum c_j \chi_j (r_j\log r_j)^{-2} - \Delta_e\log(1+B) \in r^{-\nu-2}\calC^{0,\alpha}_\eucl (M)$. 
In this reduction, it is important to note that both
\[
|\nabla^g \log(1+B)|^2_g \leq C \qquad \mbox{and} \qquad |\nabla^g \log \log r|^2_g \leq C,
\]
so we must only ensure that $h$ also has a bounded gradient. 

The Fredholm properties of the Laplacian on asymptotically Euclidean surfaces is well-known. To state
it properly, define the $2\ell$-dimensional space
\[
\calE = \{\sum_{j=1}^\ell \chi_j(\alpha_j \log r + \beta_j): \alpha_j, \beta_j \in \RR\}
\] 

\begin{lemma}
\label{lem-stand}
For any $0 < \nu < 1$, the mapping
\[
\Delta_e: \calE \oplus r^{-\nu}\calC^{2,\alpha}_\eucl(M) \longrightarrow r^{-\nu-2}\calC^{0,\alpha}_\eucl (M)
\]
is surjective and has a nullspace $K$ with  $\dim K = \ell$. 
\end{lemma}
One proof is obtained by combining Theorems 4.20, 4.26 and 7.14 in \cite{M}, but that route requires a lot 
of machinery. The $L^2$ and $L^p$ Fredholm theory is contained in \cite{LM}, see also \cite{Me}.
However, in this low-dimensional setting, it is also elementary to prove the result in H\"older spaces (including 
the existence of the supplementary space $\calE$) by direct arguments using barriers and
sequences of solutions on compact exhaustions of $M$, see \cite{Pac} .

By any of these approaches, we see that there exists a function $\tilde{h} \in \calE \oplus r^{-\nu}
\calC^{2,\alpha}_\eucl(M)$ which satisfies $\Delta_e \tilde{h} = \tilde{q}$ and 
\[
\tilde{h} = \sum_{j=1}^\ell \chi_j (\alpha_j \log r_j + \beta_j) + \eta, \qquad \eta \in r^{-\nu}\calC^{2,\alpha}_\eucl (M).
\]
This is still not enough to prove the theorem, however, because unlike the correction term $\log \log r$ used
earlier, 
\[ 
\left|\nabla^g \log r\right|^2_{g}  \approx (\log r)^2
\]
on each end. In other words, the solution $h$ of (\ref{eq:potfcn}) using this $\tilde{h}$ will not satisfy 
the conclusion of (\ref{pr:potential}) unless all the coefficients $\alpha_j$ vanish. 

Fortunately, there is flexibility in choosing the solution $\tilde{h}$ since we can add to it any element of the 
$\ell$-dimensional nullspace $K$. We now examine whether it is possible to do this in such a way that all of 
the log terms are cancelled. For this, we must study the structure of $K$ more closely. 
By definition, if $k \in K$, then
$k =  P_j \log r_j + Q_j + \calO(r^{-\nu})$ on each end $E_j$, for some constants $P_j, Q_j$. If there exists 
some $k \in K$ with $P_j = \alpha_j$, $j = 1, \ldots, \ell$, then $h = \tilde{h} - k - C$ is the solution 
we seek (the constant $C$ is chosen so that $\int f = 0$). 

\begin{lemma}
The natural map $\Phi$ defined by
\[
K \ni k \longrightarrow (P_1,Q_1, \ldots, P_\ell,Q_\ell) \in \RR^{2\ell}
\]
is injective and has image $V$ an $\ell$-dimensional subspace of $\RR^{2\ell}$ which is Lagrangian with respect to
the symplectic form
\[
\omega((P,Q), (\tilde{P},\tilde{Q})) = \sum_{j=1}^\ell (P_j \tilde{Q}_j - \tilde{P}_j Q_j).
\]
\label{le:lemma8}
\end{lemma}
\begin{proof}
The fact that $V$ is isotropic with respect to $\omega$ follows from Green's formula: if $k, k' \in K$, 
with $\Phi(k) = (P,Q)$, $\Phi(k') = (P',Q')$, then
\[
0 = \lim_{R_j \to \infty} \int_{r_j \leq R_j}  \left((\Delta_e k) k' - k (\Delta_e k')\right)\, dA_e = 
2\pi \omega((P,Q), (P',Q')).
\]
(Note that the coefficient of $\log r_j$ on each end cancels out.) But $\dim V = \ell$, so $V$ is Lagrangian.
\end{proof}

We use this as follows. Since
\[
0 = \int_M (R - \olR)\, dA = \int_M \kappa(2/\sigma) A (r \log r)^{-2}(R - \olR)\, dA_e = \int_M q\, dA_e
\]
then another application of Green's theorem gives
\[
0 = \int_M \Delta_e \tilde{h}\, dA_e = 2\pi \sum_{j=1}^\ell \alpha_j.
\]
Since $1 \in K$, its image $\Phi(1) = (0,1,0,1,\ldots,0,1)$ lies in the subspace $V$. Choose
any $(\ell-1)$-dimensional complement $K'$ to the element $1$ in $K$ and let $V' = \Phi(K')$. Then for any 
$k' \in K'$, writing $\Phi(k') = (P',Q') \in W$, we obtain as before that 
\[
\int_M \Delta_e k'\, dA_e = 2\pi \sum_{j=1}^\ell P_j'.
\]
However, by Lemma~\ref{le:lemma8}, for any such $k'$, 
\[
\omega((P',Q'), \Phi(1)) = \sum_{j=1}^\ell P_j' = 0. 
\]
Since $\Phi$ is injective, we conclude that 
\[
K' \ni k' \mapsto \Phi(k') \mapsto \{P' \in \RR^\ell: \sum_{j=1}^\ell P_j' = 0\}
\]
is bijective. This proves that we can choose the required element $k \in K$ which has
asymptotic logarithmic coefficients at all ends matching those of $\tilde{h}$, hence 
$\tilde{h} - k$ is bounded. Clearly it also has vanishing integral and bounded gradient. 

It remains to weaken the assumption on the decay properties of $A-1 = B$. Thus let $g = 
\kappa(2/\sigma) A g_c$ where $A - 1 = B \in s^{-\mu}\calC^{2,\alpha}_\hc(M)$ for some $\mu > 1$. 
We assume that $g_c$ is extended as a smooth metric on all of $M$ which is exactly hyperbolic 
on each end. Then $\Delta_g f = R - \olR$ is equivalent to 
\[
\Delta_c f = \kappa(2/\sigma)A(R - \olR) = \kappa(2/\sigma)A \left( (R + \kappa(\sigma)) - 
(\kappa(\sigma) + \olR) \right);
\]
by the computations in \S 2, modulo $s^{-\mu}\calC^{0,\alpha}_\hc$, the function $A$ can
be replaced by $1$ and $R + \kappa(\sigma) \equiv 0$, so the right hand side has the form 
\[
- \kappa(2 + 2\olR/\sigma) + \zeta, \qquad \zeta \in s^{-\mu}\calC^{0,\alpha}_\hc.
\]
Notice also that
\[
\int \kappa(2/\sigma) A (R - \olR)\, dA_c = \int (R - \olR)\, dA = 0.
\] 

Now consider $h_j :=\chi_j (2/\sigma_j)A(R -\olR)$ as a function on the simple quotient $X = {\mathbb H}^2/\ZZ 
\cong \RR \times S^1$ with hyperbolic metric $ds^2 + e^{-2s}d\theta^2$ (which we also call $g_c$). 
For convenience, drop the index $j$ for the moment. To produce a solution $w$ to $\Delta_c w = h$ 
on $X$ which satisfies $|w| \leq C s^{-\mu + 1}$ for $s \geq 0$, change the metric conformally one 
more time. Let $\tau = e^s$, so that $\tau^2 g_c = g_{\mathrm{cyl}} = d\tau^2 + d\theta^2$ 
on the half-cylinder $\RR^+ \times S^1$. Then $\Delta_c w = h$ is equivalent to $\Delta_{\mathrm{cyl}} w = \tau^{-2}h$.
Note that $|\tau^{-2}h| \leq C \tau^{-2}(\log \tau)^{-\mu}$. Decompose $w = w_0 + w^\perp$ where $w_0$ 
is the Fourier mode of order zero on the $S^1$ factor and $w^\perp$ is the sum of all other Fourier components, 
and similarly for $\tau^{-2}h$.  It is easy to find $w_0$ by integrating in from $\tau = \infty$, and to 
check that $|w_0| \leq C (\log \tau)^{-\mu + 1}$ for $\tau \geq 2$. On the other hand, we claim that the
unique $L^2$ solution $w^\perp$ satisfies $|w^\perp| \leq C \tau^{-1}$ for $\tau \geq 1$ (this is a very 
crude bound, but is sufficient for what we need). To see this, take the Fourier transform in $\tau$ and 
expand in Fourier series in the $S^1$ factor. Denote the corresponding Fourier coefficients
of $w^\perp$ and $\tau^{-2}h$ by $\hat{w}(\lambda)$ and $\hat{k}_j(\lambda)$, respectively. 
By the assumed regularity of $h$ and the decay of $\tau^{-2}h$, we certainly have a uniform bound
$|\hat{k}_j(\lambda)| \leq C$ independent of $j$ and $\lambda$. Hence 
\[
|\hat{w}_j(\lambda)|\leq C/(j^2 + \lambda^2) \qquad \mbox{and} \qquad
|\del_\lambda \hat{w}_j(\lambda)| \leq C/(j^2 + \lambda^2)^2.
\] 
Now write $i\tau e^{i\tau \lambda} = \del_\lambda e^{i\tau \lambda}$ and integrate by parts in $\lambda$ to obtain
\[
|\tau w^\perp(\tau,\theta)| \leq \sum_{j \neq 0} \int_{-\infty}^\infty \frac{C}{(j^2 + \lambda^2)^2}\, d\lambda
\leq C'
\]
as claimed. 

Finally, reclaiming the index $j$, transplant $\chi_j w_j$ back to $E_j \subset M$ and define 
\[
h_1 = \kappa(2/\sigma) A (R - \olR) - \sum_j \Delta_c (\chi_j w_j). 
\]
This is compactly supported, and by Green's formula, $\int \Delta_c (\chi_j w_j)\, dA_c = 0$,
so $\int_M h_1 \, dA_c = 0$. Now apply the first part of the proof to find a function $v$ such that 
$\Delta_c v = h_1$, $|v| \leq C r^{-\nu}$ on each end, and such that
$\int v = 0$, $|\nabla v| \leq C$. The solution to the original problem is then 
$f = v + \sum \chi_j w_j$.

This proves Proposition \ref{pr:potential} in generality.   \hfill  $\Box$

\section{Long-time existence and convergence}
The remaining tasks are to prove that the normalized Ricci flow applied to the metric $g_0$
exists for all time and converges to a hyperbolic metric. 

We shall use the maximum principle in the following form.
\begin{lemma}
\label{lem-max}
Let $g(t)$ be a family of complete Riemannian metric on a noncompact manifold $M$, defined for $t\in [0,T)$, which
varies smoothly in $t$, and such that each $g(t)$ is quasi-isometric to $g(0)$. Let $f(\cdot,t)$ be a smooth bounded 
function on $M\times [0,T)$ such that $f(x,0) \ge 0$, and
\[
\frac{\del\,}{\del t} f = \Delta_{g(t)} f + Q(f,x,t), \qquad \mbox{where}\qquad 
Q(f,x,t) \ge 0\quad \mbox{when}\quad f \le 0.
\]
Then $f(x,t) \ge 0$ on $M\times [0,T)$.
\end{lemma}
The proof can be found in \cite{Sh}, \cite{Sh1}.

We now take advantage of the fact that in dimension $2$, every metric is K\"ahler, so we can regard
(\ref{equation-RF}) as the K\"ahler-Ricci flow; in other words, as in \cite{Cha}, using the complex structure associated
to this conformal class, we write
\begin{equation}
g_{i\bar{\jmath}} = (g_0)_{i\bar{\jmath}} + \del_i \del_{\bar{\jmath}} \phi
\label{eq:gg0}
\end{equation}
where $\phi(t)$ is a solution to the complex parabolic Monge-Ampere equation
\begin{equation}
\label{eq-metric-pot}
\frac{\del\,}{\del t}\phi = \log \frac{\det((g_0)_{k\bar{\ell}} + \phi_{k\bar{\ell}}))}
{\det((g_0)_{k\bar{\ell}})} + \olR \, \phi - 2f_0; \qquad \phi(x,0) = 0.
\end{equation} 
here $f_0$ is the solution of  (\ref{eq:potfcn}) as provided by Proposition \ref{pr:potential}.

Now define $f = -\frac12 \del_t \phi$. Differentiating (\ref{eq-metric-pot}) with respect to $t$ gives
\begin{equation}
\label{eq-ev-pot}
\frac{\del\, }{\del t} f = \Delta f + \olR\, f,  \qquad \left. f \right|_{t=0} = f_0.
\end{equation} 
We claim that, up to a summand depending only on $t$, $f$ is the potential function for $g$ for all 
values of $t$, i.e.\ $\Delta_{g(t)} f(t,\cdot) = R_t - \olR_t$ and $|\nabla_t f| \leq C$. To see this,
first recall that 
\[
R \, g_{i\bar{\jmath}} = R_{i\bar{\jmath}} = -\partial_i\partial_{\bar{\jmath}}\log\det g_{k\bar{\ell}} \, ;
\]
next apply $\del_i\del_{\bar{\jmath}}$ to  (\ref{eq-metric-pot}) to get
\begin{eqnarray*}
-2 \del_i\del_{\bar{\jmath}}f &=& - Rg_{i\bar{\jmath}} + R_0 (g_0)_{i\bar{\jmath}} 
+ \olR \del_i\del_{\bar{\jmath}}\phi - 2\del_i\del_{\bar{\jmath}} f_0 \\
&=& - (R - \olR)g_{i\bar{\jmath}} + (R_0 - \olR_0)(g_0)_{i\bar{\jmath}} - 2\del_i\del_{\bar{\jmath}}f_0,
\end{eqnarray*}
using (\ref{eq:gg0}) and the fact that the average value $\olR = \olR_0$ is independent of $t$. 
Taking the trace, and reverting back to conformal notation, $g = u g_0$, yields
\begin{equation}
\label{eq-trace}
\Delta f = R - \olR - \frac{1}{u} ( (R_0 - \olR_0) - \Delta_0 f_0) = R - \olR.
\end{equation}

Now, following \cite{Ha}, define the function $h := \Delta f + |\nabla f|^2$ and the symmetric $2$-tensor 
$Z := \nabla^2 f - \frac{1}{2}\Delta f\cdot g$. A simple computation shows that $h$ satisfies
\[
\frac{\del h}{\del t} = \Delta h - 2|Z|^2 - h.
\]
In order to apply the maximum principle, we must show that $|h|$ is bounded on $M$ for each time $t$.
This depends, in turn, on 
\begin{lemma}
If the flow exists for $0 \leq t < T$, then $\sup_M |\nabla f|$ is finite for every $t < T$. 
\label{le:nablaf}
\end{lemma}
\begin{proof}
We compute the evolution equation for $|\nabla f|^2$: since
\[
\frac{\del\, }{\del t}g^{ij} = -g^{ip}g^{jq}\frac{\del\, }{\del t}g_{pq} = (R - \olR) g^{ij},
\]
we have
\begin{equation}
\label{eq-nabla-f}
\begin{split}
\frac{\del\, }{\del t}|\nabla f|^2 &= \frac{\del\, }{\del t}g^{ij}\nabla_if\nabla_j f + 2g^{ij}
\nabla_i\frac{\del f}{\del t}\nabla_j f \\
&= (R+1)|\nabla f|^2 + 2g^{ij}\nabla_i(\Delta f - f)\nabla_j f \\
&= (R+1)|\nabla f|^2 +2\langle\nabla R, \nabla f\rangle - 2|\nabla f|^2 \\
&\le C_1|\nabla f|^2 + C_2t^{-1/2}|\nabla f|,
\end{split}
\end{equation}
where the last inequality uses Shi's estimates $|R(x,t)| \le C(t_1)$ and $|\nabla R| \le C(t_1)t^{-1/2}$ if the
flow exists on $M\times (0,t_1]$. For each $x$, fix $\tau = \tau(x)$ small, and denote $D(x,\tau) = 
\sup_{[0,\tau]} |\nabla f(x,t)|^2$; thus $\sup_M D(x,0) < \infty$. Integrating (\ref{eq-nabla-f}) from $0$ to $\tau$ gives
\[
\begin{split}
D(x,\tau) & \leq C_0 + C_1 \tau D(x,\tau) + C_2 D(x,\tau) \int_0^{\tau} t^{-1/2}\, dt \\
& \leq C\left(1 + \tau D(x,\tau) + \sqrt{\tau D(x,\tau)}\right) \leq C(1+\sqrt{\tau} D(x,\tau)),
\end{split}
\]
assuming that $\tau < 1$ and $D(x,\tau) > 1$. The constant $C$ is independent of $x$, hence so is
$\tau = (4C^2)^{-1}$, and with this $\tau$ we obtain the uniform upper bound
\begin{equation}
\label{eq-first}
\sup_{M \times [0,\tau]} |\nabla f(x,t)|^2 \leq C.
\end{equation}

Finally, for $\tau \leq t < T$, 
\[
\frac{\del\, }{\del t}|\nabla f|^2 \le C_1|\nabla f|^2 +C_2 \tau^{-1/2}|\nabla f|
\]
so integrating from $\tau$ to any other value $t < t_0$ and using (\ref{eq-first}) gives
\[
\sup_M |\nabla f(\cdot,t)|^2 \le C(t),
\]
as desired, so long as $t < T$. 
\end{proof}

Since $|R(x,t)|$ is bounded in $x$ for $t < T$, we also have
\begin{equation}
\label{eq-delta-f}
\Delta f (x,t) = R - \olR \le C(t) \qquad \mbox{for} \quad t < T,
\end{equation}
and hence 
\[
\sup_M h(\cdot,t) \le C(t) < \infty \qquad \mbox{for}\quad t < T
\]
as well. 

It now follows directly from the equation for $h$ and the maximum principle that
\[
h \le Ce^{-t}
\]
for all $t < T$, where now the constant $C$ is independent of $T$. Applying the maximum principle to (\ref{eq:eqforR}) 
gives a $T$-independent lower bound on $R$, so altogether, since $R = h - |\nabla_g f|^2 - 1$, we have
\begin{equation}
\label{equation-R-neg}
-C \le R \le Ce^{-t} + \olR,
\end{equation}
where $C$ is independent of $T$. 

We have now proved that if the flow exists on $[0,T)$, then there is an a priori estimate for
$R$ on this interval, and it is standard that this implies that the maximal interval
of existence is $(0,\infty)$. Indeed, integrating the flow equation itself gives exponential
upper and lower bounds for the metric. Taking repeated covariant derivatives of (\ref{eq:eqforR})
and bootstrapping gives a priori estimates for all higher derivatives of $R$ as well, and
from this it follows that $g(t)$ exists for all $t < \infty$. 
Now, just as in the compact case \cite{Ha}, (\ref{equation-R-neg}) shows that
$R < 0$ when $t$ is greater than some $\overline{t}$, so a further application of the 
maximum principle to (\ref{eq:eqforR}) yields 
\[
-e^{-\e (t-\overline{t})} \le \olR - R \le Ce^{-t}
\]
when $t \geq \overline{t}$. Hence $R \to \olR$ exponentially. Integrating the flow equation one
final time, from any $\tilde{t} \geq \overline{t}$ to infinity, shows that $g(t)$ converges 
to the uniformizing metric of constant Gauss curvature $2\pi \chi(M)/A_0$.

\end{document}